\documentclass[a4paper,10pt]{amsart}

\usepackage{amsaddr} 

\usepackage[utf8]{inputenc}

\usepackage[backend=biber, style=alphabetic, sorting=ynt]{biblatex}
\bibliography{../paper.bib}

\usepackage{microtype}
\usepackage{xcolor}
\usepackage{hyperref}
\usepackage{listings}

\hypersetup{ 
	pdftitle={pySymmPol},
	pdfsubject={Computational Physics,
      Mathematical-Physics, High Energy Physics}, 
	pdfauthor={Thiago Araujo},
	pdfkeywords={python},
	colorlinks=true, 
    linkcolor=red, 
    citecolor=red, 
    urlcolor=violet,
    linktocpage=true
  }

\definecolor{myPurple}{RGB}{90, 74, 120}
\definecolor{myBlue}{RGB}{15, 75, 110}
\definecolor{myRed}{RGB}{191,97,106}

\begin{document}

\title[PySymmPol]{PySymmPol: Symmetric Polynomials in Python (v0.1.2)}

\author{Thiago Araujo}

\address{\noindent 
Instituto de Física Teórica, UNESP-Universidade Estadual Paulista,
R. Dr. Bento T. Ferraz 271, Bl. II, Sao Paulo 01140-070, SP, Brazil\\
\&
Instituto de Física, Universidade de S\~ao Paulo,
Rua do Matão Travessa 1371, 05508-090 São Paulo, SP. Brazil
}

\email{\href{tr.araujo@unesp.br}{tr.araujo@unesp.br}}

\keywords{python, symmetric polynomials}
\date{\today}

\maketitle

\setcounter{tocdepth}{1}
\tableofcontents

\section{Summary}

\textbf{PySymmPol} is a Python package designed for efficient manipulation of
symmetric polynomials. It provides functionalities for working with
various types of symmetric polynomials, including elementary,
homogeneous, monomial symmetric, (skew-) Schur, and Hall-Littlewood
polynomials. In addition to polynomial operations, \textbf{PySymmPol} offers
tools to explore key properties of integer partitions and Young
diagrams, such as transposition, Frobenius coordinates, characters of
symmetric groups and others. 

This package originated from research conducted in the field of
integrable systems applied to string theory, and the AdS/CFT (Anti-de
Sitter/Conformal Field Theory) correspondence. \textbf{PySymmPol} aims to
facilitate computational tasks related to symmetric polynomials and
their applications in diverse fields.

Instructions for instalation and tutorials for the different modules
are in the documentation page. 
\begin{itemize}
    \item Github repository: \href{https://github.com/thraraujo/pysymmpol}{github.com/thraraujo/pysymmpol}
    \item Documentation: \href{https://thraraujo.github.io/pysymmpol}{thraraujo.github.io/pysymmpol}
    \item PyPI: \href{https://pypi.org/project/pysymmpol}{pypi.org/project/pysymmpol}
\end{itemize}
Instalation with pip 
\begin{lstlisting}[language=Python]
  $ pip install pysymmpol
\end{lstlisting}

\section{Statement of need}

Symmetric polynomials play a crucial role across various domains of
mathematics and theoretical physics due to their rich structure and
broad applications. They arise naturally in combinatorics \cite{Macdonald:1998}, 
representation theory \cite{Fulton:2004}, algebraic geometry, and mathematical
physics \cite{Babelon:2003, Wheeler:2010}. These polynomials encode
essential information about symmetries and patterns, making them
indispensable in the study of symmetric functions and their
connections to diverse mathematical structures. Moreover, symmetric
polynomials find extensive applications in theoretical physics,
particularly in quantum mechanics, statistical mechanics, and quantum
field theory. Their utility extends to areas such as algebraic
combinatorics, where they serve as powerful tools for solving
combinatorial problems and understanding intricate relationships
between different mathematical objects. Thus, tools and libraries like
PySymmPol provide researchers and practitioners with efficient means
to explore and manipulate symmetric polynomials, facilitating
advancements in both theoretical studies and practical applications.
 
Let's now provide a brief overview of the addressed problems. In the
next section, we will explain how the package can be utilized.

\subsection{Integer Partitions and Young Diagrams}

For integer partitions, we employ two primary representations. The
first is the conventional representation as a monotonic decreasing
sequence, denoted as $\lambda = (\lambda_1, \dots, \lambda_m)$, where
$\lambda_i \geq \lambda_{i+1}$. This representation enables basic
manipulations of integer partitions, including determining the number
of rows, columns, boxes, diagonal boxes, and Frobenius coordinates
within the associated Young Diagram.

Additionally, Young diagrams are represented as cycles within the
symmetric group $\mathfrak{S}_N$, denoted as $\lambda = (1^{k_1}
2^{k_2} \dots)$, where $N = \sum_i\lambda_i = \sum_j j k_k$. These
cycles serve as representatives of the conjugacy class of
$\mathfrak{S}_N$, and as such, we refer to them as conjugacy class
vectors. This alternative representation offers insights into the
structural properties of integer partitions and facilitates their
manipulation within the framework of group theory.

Our keen interest in exploring these representations stems from their
relevance in the realm of two-dimensional Integrable and Conformal
Field Theories, see \cite{Babelon:2003, Marino:2005, Okounkov:2006} for
a small sample of physical problems that require tools to manipulate
these (and others) combinatorics objects. In this context, bosonic
states are labeled by conjugacy class vectors, highlighting the
importance of understanding and manipulating such
representations. Conversely, fermions are represented using the
standard representation, emphasizing the need to bridge the conceptual
gap between these distinct frameworks. Investigating these
representations not only sheds light on the fundamental properties of
fermion states within CFTs but also provides valuable insights for
theoretical developments in quantum field theory and related fields.

\subsection{Symmetric Polynomials}

One of the main goals of the package is to provide the definitions of the 
polynomials in terms of the Miwa coordinates, or power sums, 
$$ t_j = \frac{1}{j} \sum_{i=1}^N x_i^j $$ where $\vec{x} = 
(x_1, \dots, x_N)$ and $\mathbf{t} = (t_1, t_2, \dots)$.

In terms of these coordiantes, the \emph{Complete Homogeneous} $h_n(\mathbf{t})$ 
and \emph{Elementary Symmetric Polynomials} $e_n(\mathbf{t}) = 
(-1)^n h_n(-\mathbf{t})$ are defined via 
$$ h_n(\mathbf{t}) = \sum_{k_1 + 2k_2+ \cdots = n} 
\frac{t_1^{k_1}}{k_1}\frac{t_2^{k_2}}{k_2} \cdots $$
From these expressions, we obtain the \emph{(skew-) Schur polynomials}
$s_{\lambda/\mu}(\mathbf{t})$, where 
$\lambda$ and $\mu$ are integer partitions, via \emph{Jacobi-Trudi identity}
$$ s_{\lambda/\mu} = \det_{p,q}(h_{\lambda_q - \mu_p - q + p}(\mathbf{t})) \; . $$

The \emph{Hall-Littlewood polynomials} are defined via
$$
P_{\lambda}(x_1, \dots, x_N; Q) = \prod_{i\geq 0} \prod_{j=1}^{p(i)}
\frac{1-Q}{1-Q^j} \sum_{\omega \in \mathfrak{S}_N} \omega\left(
x_1^{\lambda_1}\cdots x_n^{\lambda_n} \prod_{i<j} \frac{x_i - Q
  x_j}{x_i - x_j} \right)
$$
where $p(i)$ is the number of rows of size $i$ in $\lambda$, and
$\mathfrak{S}_N$ is the symmetric group. The limit $Q=0$ in the
Hall-Littlewood polynomials gives the Schur polynomials, while $Q=1$
returns the \emph{Monomial Symmetric Polynomials} $m_\lambda(x_1,
\dots, x_N)$. For more information on these topics, see the references
below.

\section{Basic Usage}

The \textbf{PySymmPol} package has seven main classes for manipulating
various symmetric polynomials.

\subsection{YoungDiagram} and ConjugacyClass

For the construction and manipulation of Young diagrams, we need to import 
the YoungDiagram and the ConjugacyClass classes. 
\begin{lstlisting}[language=Python]
from pysymmpol import YoungDiagram, ConjugacyClass
\end{lstlisting}

The distinction between these two classes lies in the representations of the diagrams 
they handle. The YoungDiagram class represents diagrams using a monotonic decreasing sequence, 
while the ConjugacyClass class represents them as a sequence representing the cycle 
of the symmetric group $\mathfrak{S}_N$. For example, let us consider the partition
$(3,2,1)$, which is represented as a tuple in the YoungDiagram class and 
as a dictionary in the ConjugacyClass class $\{1: 1, 2: 1, 3: 1\}$, respectively.
\begin{lstlisting}[language=Python]
young = YoungDiagram((3,2,1))
conjugacy = ConjugacyClass({1: 1, 2: 1, 3: 1})
\end{lstlisting}
Both objects describe the same mathematical entity, the partition \(6=3+2+1\).
In fact, 
we have the usual pictorial representation 
\begin{lstlisting}[language=Python]
young.draw_diagram(4)

conjugacy.draw_diagram(4)
\end{lstlisting}
that give the same output (the argument 4 means that we draw the octothorpe, and
there are 4 other symbols available).
\begin{lstlisting}[language=Python]
#
# #
# # #

#
# #
# # #
\end{lstlisting}
Further details on the other functionalities can be found in the tutorial.
The \textrm{ACCELASC} algorithm~\cite{Kelleher:2009} greatly improved the
speed of these methods.

\subsection{Homogeneous and Elementary Polynomials}

These classes can be initialized as 
\begin{lstlisting}[language=Python]
from pysymmpol import HomogeneousPolynomial, ElementaryPolynomial
from pysymmpol.utils import create_miwa
\end{lstlisting}
We also imported the function \verb|create_miwa| from the \verb|utils| module for convenience. 
Now, let us create the polynomials at level $n=3$. We can instantiate 
them and find their explicit expressions using the \verb|explicit(t)| method, where \verb|t| 
represents the Miwa coordinates. The block
\begin{lstlisting}[language=Python]
n = 3
t = create_miwa(n)

homogeneous = HomogeneousPolynomial(n)
elementary = ElementaryPolynomial(n)
print(f"homogeneous: {homogeneous.explicit(t)}")
print(f"elementary: {elementary.explicit(t)})
\end{lstlisting}
gives the output 
\begin{lstlisting}[language=Python]
homogeneous: t1**3/6 + t1*t2 + t3
elementary: t1**3/6 - t1*t2 + t3
\end{lstlisting}

\subsection{Schur Polynomials}

To create Schur polynomials, we first need to instantiate a partition 
before defining the polynomial itself. Let us use the Young diagram we 
considered a few lines above, then
\begin{lstlisting}[language=Python]
from pysymmpol import YoungDiagram
from pysymmpol import SchurPolynomial
from pysymmpol.utils import create_miwa
\end{lstlisting}
The YoungDiagram class includes a getter method for retrieving the number 
of boxes in the diagram, which we utilize to construct the Miwa coordinates. 
Subsequently, the SchurPolynomial class is instantiated using the Young diagram. Then
\begin{lstlisting}[language=Python]
young = YoungDiagram((3,2,1))
t = create_miwa(young.boxes)

schur = SchurPolynomial(young)

print(f"schur: {schur.explicit(t)}")
\end{lstlisting}
gives
\begin{lstlisting}[language=Python]
schur: t1**6/45 - t1**3*t3/3 + t1*t5 - t3**2
\end{lstlisting}
The documentation and tutorial contain examples demonstrating how to find 
\textbf{skew-Schur polynomials.}

\subsection{Monomial Symmetric Polynomials}

For Monomial symmetric polynomials, we have a similar structure. 
\begin{lstlisting}[language=Python]
from pysymmpol import YoungDiagram
from pysymmpol import MonomialPolynomial
from pysymmpol.utils import create_x_coord
\end{lstlisting}
The only difference is the function \verb|create_x_coord| from the
\verb|utils| module. Therefore,
\begin{lstlisting}[language=Python]
young = YoungDiagram((3,2,1))

n = 3
x = create_x_coord(n)

monomial = MonomialPolynomial(young)

print(f"monomial: {monomial.explicit(x)}")
\end{lstlisting}
gives the output 
\begin{lstlisting}[language=Python]
monomial: x1*x2*x3*(x1**2*x2 + x1**2*x3 + x1*x2**2 + x1*x3**2 + x2**2*x3 + x2*x3**2)
\end{lstlisting}

\subsection{Hall-Littlewood Polynomials}

In addition to partitions, for the Hall-Littlewood polynomials, 
we also require the deformation parameter \(Q\) (as \(t\) has been used to 
denote the Miwa coordinates).
\begin{lstlisting}[language=Python]
from sympy import Symbol
from pysymmpol import YoungDiagram
from pysymmpol import HallLittlewoodPolynomial
from pysymmpol.utils import create_x_coord
\end{lstlisting}
The method \verb|explicit(x, Q)| needs another argument. Finally, the code
\begin{lstlisting}[language=Python]
Q = Symbol('Q')
young = YoungDiagram((3,2,1))

n = 3
x = create_x_coord(n)

hall_littlewood = HallLittlewoodPolynomial(young)

print(f"hall-littlewood: {hall_littlewood.explicit(x, Q)}")
\end{lstlisting}
gives
\begin{lstlisting}[language=Python]
hall-littlewood: x1*x2*x3*(-Q**2*x1*x2*x3 - Q*x1*x2*x3 + x1**2*x2 + x1**2*x3 + x1*x2**2 + 2*x1*x2*x3 + x1*x3**2 + x2**2*x3 + x2*x3**2)
\end{lstlisting}

\section{State of the field and target audience}

To the best of the author's knowledge, there are few open-source
software solutions dedicated to similar problems, with SageMath
\cite{sagemath} standing out as a notable exception. While there are
overlaps between the problems addressed and some implementations
available in SageMath, significant differences exist. One prominent
dissimilarity is that SageMath constitutes a collection of open-source
mathematical software, whereas \textbf{PySymmPol} is a Python package. A
notable characteristic of \textbf{PySymmPol} is its utilization of power
sums and its Pythonic nature, which minimizes dependencies.

Furthermore, implementations in SageMath primarily emphasize
applications in combinatorics. Although the author extensively used
SageMath, the development of \textbf{PySymmPol} stemmed from the need for a
more physics-oriented software. For instance, it offers features such
as the translation of standard notation of partitions (representing
fermionic states in certain 2D field theories) and the conjugacy class
notation of partitions (representing bosonic states).

Another significant distinction between our implementation and
SageMath is the use of Miwa coordinates (or power sums) in the
definitions. This aspect proves advantageous for physicists and
mathematicians involved in statistical physics, quantum field theory,
and integrable systems.

\section{Acknowledgements}

I would like to thank Fapesp for financial support, grant
\textbf{2022/06599-0}.  I would also like to thank the Free and Open Source
Software community.

\printbibliography
     
\end{document}